\documentclass[twoside]{article}

\usepackage{tikz}\usetikzlibrary{hobby}

\usepackage{problems}

\usepackage[colorlinks=true,
citecolor=black,
linkcolor=black,
anchorcolor=black,
filecolor=black,
menucolor=black,
urlcolor=black,
pdftitle={A gallery of open problems in geometry that keep me puzzled and amused},
pdfsubject={Geometry},
pdfauthor={Anton Petrunin}
]{hyperref}

\makeindex

\begin{document}

\title{
A gallery of open problems in geometry\\
that keep me puzzled and amused
}
\author{Anton Petrunin}
\date{}
\maketitle

This is a collection of open questions in geometry that I think of as \textit{puzzles}:
they stick to my brain --- I see grips, but no spare hands.
Puzzle-charm is the only criterion for inclusion here;
importance is ignored.

These puzzles grew out of my collection of solved puzzles~\cite{petrunin2021pigtikal}.
I kept it aside since things around open questions change rapidly.

Please send me an email if you solve one of the puzzles, or see a good approach.

\renewcommand\thesubsection{\arabic{subsection}}

\section*{Homemade}

\subsection{Two-convexity}
\label{Two-convexity}

Note that an open connected set $\Omega\subset\mathbb R^n$ is convex (in the usual sense) if the following holds:
for any triangle $\triangle$, if two sides of $\triangle$ lie in $\Omega$, then $\triangle\subset\Omega$.
This observation motivates the definition of two-convexity used in the following problem.

\begin{pr}
Let $\Omega$ be a simply connected open set in $\mathbb R^n$.
Assume that $\Omega$ is two-convex;
that is, if 3 faces of a 3-simplex $\triangle$ lie in $\Omega$, then $\triangle\subset\Omega$.

Is it true that any component of the intersection of $\Omega$ with any 2-plane is simply connected?
\end{pr}

If the boundary of $\Omega$ is a smooth hypersurface, then the answer is ``yes''.
Also, the answer is ``yes'' if $n=3$.

To prove the first statement, note that two-convexity implies that
at most one of the principal curvatures of the boundary is negative.
Then follow the Morse-type argument in Lefschetz theorem \cite[see Section~$\tfrac12$ in][]{gromov-SGMC}.
In the case $n=3$, one can mimic the Morse-type argument;
it also follows from \cite[7.17 in][]{AKP-invitation-2}.

The first argument also proves that \textit{if $\Omega$ is a simply connected open two-convex set in $\mathbb R^n$ with a smooth boundary, then any component of the intersection of $\Omega$ with any affine subspace (not necessarily 2-dimensional) is simply connected.}
The example below shows that a general 2-convex open set does not have this property.

One may hope to prove the statement by constructing an approximation of $\Omega$ by two-convex domains with smooth boundaries.
Again, the example below shows that such an approximation does not exist for $n\ge 4$.

\parit{Example.}
We will construct a two-convex simply connected open set $\Omega$ in $\mathbb R^4$ such that the intersection of $\Omega$ with some hyperplane is not simply connected.

Set
$$\Pi=\{\,(x,y,z,t)\in\mathbb R^4\mid\,y< x^2\}.$$
Let $\Pi'$ be the image of $\Pi$ under a generic rotation of $\mathbb R^4$,
say $(x,y,z,t)\z\mapsto(z,t,x,y)$.

Since $\Pi$ is open and two-convex, so is the intersection $\Omega=\Pi\cap \Pi'$.
Note that we can choose $(x,y,z,t)$-coordinates so that $\Omega$ is an epigraph for a function $f\colon\mathbb R^3\to\mathbb R$ like
$$f=\max\{\alpha_1-\beta_1^2,\alpha_2-\beta_2^2 \},$$
where $\alpha_i$ and $\beta_i$ are linear functions on the $(x,y,z)$-coordinate subspace.
Thus, in particular, $\Omega$ is contractible.

Let $L_{t_0}$ be the intersection of $\Omega$ with hyperplane $t=t_0$;
it is a complement of two convex parabolic cylinders in general position.
If these cylinders have a point of intersection, then $\pi_1 L_{t_0}\z=\mathbb Z$.

The discussion above implies that $f$ cannot be approximated by smooth functions such that their Hessians have at most one negative eigenvalue at all points.

\parit{Remark.}
Two-convexity shows up in comparison geometry --- the maximal open flat sets in the manifolds of nonnegative or nonpositive curvature are two-convex \cite{panov-petrunin}.

\subsection{Braid space}

Let us define \emph{braid space} $\mathbb{B}^n$ as the universal covering of $\CC^n$ infinitely branching along the complex hyperplanes $z_i=z_j$ for $i\ne j$.
We suppose that $\CC^n$ is equipped with the canonical Euclidean metric,
and $\mathbb{B}^n$ is equipped with the lifted length metric.

In other words, if $\mathbb{W}^n$ denotes the complement of complex hyperplanes $z_i=z_j$ in $\CC^n$,
then the braid space $\mathbb{B}^n$ is the completion of the universal metric covering of $\mathbb{W}^n$.

The following problem is a special case of a conjecture formulated by Daniel Allcock~\cite[Conjecture 1.4]{allcock}, which is related to an earlier conjecture of Ruth Charney and Michael Davis \cite[Conjecture 3]{charney-davis-95}, which in turn is motivated by the $K(\pi,1)$-conjecture \cite{vanderlek,paris}.
Dmitri Panov and the author are also responsible for a more general version of Allcock's conjecture \cite[Conjecture 1.1]{panov-petrunin}.

\begin{pr}
Is it true that $\mathbb{B}^n$ is $\mathrm{CAT}(0)$ for any $n$?
\end{pr}

The answer is ``yes'' for $n\le 3$.
The case $n=2$ is trivial, and the case $n=3$ follows from the theorem of Ruth Charney and Michael Davis about branched coverings of the 3-sphere \cite{charney-davis};
another proof is given by Dmitri Panov and the author \cite{panov-petrunin-ramification}.
The case of the real braid space with 4 strings was settled by Ruth Charney thru a detailed analysis of
short galleries \cite[Theorem~5.4]{charney2004}.
The \emph{real braid space} is defined as the inverse image of the real subspace $\RR^n\subset \CC^n$ under the projection $\mathbb{B}^n\to \CC^n$.
This is a convex subset of  $\mathbb{B}^n$; therefore, the question about real braid space is a special case of our problem.

The space $\mathbb{B}^n$ has plenty of convex functions, which might help to understand its geometry.
The following construction is inspired by the argument given by Ian Agol \cite{agol-braid}:
(1) choose $n$ distinct points in $\CC$ (that is, an element of $\mathbb{W}^n$),
(2) choose a closed curve $z$ in its complement,
(3) consider the function $\mathbb{B}^n\to \RR$  that returns the greatest lower bound of lengths of closed curves that are homotopic to $z$.

Existence of convex functions forbids geodesic loops in $\mathbb{B}^n$; this would also follow if $\mathbb{B}^n$ is $\mathrm{CAT}(0)$.
To show that $\mathbb{B}^n$ is $\mathrm{CAT}(0)$ one has to forbid geodesic digons \cite[6.6]{AKP-invitation-2}.

\subsection{Quotients of Hilbert space}

\begin{pr} Suppose $R$ is a compact simply connected Riemannian manifold that is isometric to a quotient of the Hilbert space by a group of isometries (or, more generally, $R$ is the target of a Riemannian submersion from the Hilbert space).
Is it true that $R$ is isometric to a biquotient?
That is, is it true that $R$ is a quotient of a compact Lie group $G$ with a bi-invariant metric by a group of isometries?

\end{pr}

Any biquotient can appear as a quotient of the Hilbert space by a group of isometries.
Let us sketch the argument suggested by Alexander Lytchak; see \cite{lebedeva-petrunin-zolotov}.
Another proof follows from the construction of Chuu-Lian Terng and Gudlaugur Thorbergsson given in \cite[Section 4]{terng-thorbergsson}.

Denote by $G^n$ the direct product of $n$ copies of $G$.
Consider the map $\phi_n\:G^n\to G/\!\!/H$ defined by
\[\phi_n\:(\alpha_1,\dots,\alpha_n)\mapsto [\alpha_1\cdots\alpha_n]_H,\]
where $[x]_H$ denotes the $H$-orbit of $x$ in $G$.

Note that $\phi_n$ is a quotient map for the action of $H\times G^{n-1}$ on $G^n$ defined by
\[(\beta_0,\dots,\beta_n)\cdot(\alpha_1,\dots,\alpha_n)
=
(\beta_0\cdot \alpha_1\cdot\beta_1^{-1},\beta_1\cdot\alpha_2\cdot\beta_2^{-1},\dots,\beta_{n-1}\cdot\alpha_n\cdot\beta_n^{-1}),\]
where $\beta_i\in G$ and $(\beta_0,\beta_n)\in H<G\times G$.

Denote by $\rho_n$ the product metric on $G^n$ rescaled by the factor $\sqrt{n}$.
Note that the quotient $(G^n,\rho_n)/(H\times G^{n-1})$ is isometric to $G/\!\!/H\z=(G,\rho_1)/\!\!/H$.
Let $\phi_n\:(G^n,\rho_n)\to G/\!\!/H$ be the corresponding quotient map;
clearly, $\phi_n$ is a submetry.

As $n\to\infty$, the curvature of $(G^n,\rho_n)$ converges to zero and its injectivity radius goes to infinity.
Therefore, the ultralimit of $(G^n,\rho_n)$, with the identity element as the basepoint, is a Hilbert space $\HH$, and the submetries $\phi_n$ ultraconverge to a submetry $\phi\:\HH\to G/\!\!/H$.

\subsection{Nested convex surfaces}

\begin{pr}
Describe the Riemannian metrics on $\mathbb{S}^n$ that are isometric to smooth \emph{nested convex surfaces};
that is, a complete smooth convex hypersurface
in a complete smooth convex hypersurface
in ...
in a Euclidean space.
\end{pr}

If $n=2$, then by Alexandrov's embedding theorem these are all Riemannian metrics with nonnegative curvature.

Direct calculations show that the metric has \emph{nonnegative cosectional curvature};
the latter means that, at each point, the curvature tensor can be expressed as a linear combination with nonnegative coefficients of the curvature tensors of $\mathbb{S}^2\times \RR^{n-2}$.
(In dimension three, nonnegative cosectional curvature is equivalent
to nonnegative sectional curvature; in dimension 4, it is equivalent to nonnegative curvature operator.)

It might happen that any metric with nonnegative cosectional curvature on $\mathbb{S}^n$ is isometric to a nested convex surface.
(It's embarrassing, but I don't understand the situation with the Berger spheres.)
If one drops completeness from the definition of nested convex surfaces
and assumes that the cosectional curvature is strictly positive,
then the answer is ``yes''~\cite{petrunin-poly}.

\subsection{Maximal finite subgroups}

\begin{pr}
Suppose that a group $\Gamma$ acts effectively, isometrically, and properly discontinuously on Euclidean space $\mathbb{E}^n$.
Show that $\Gamma$ contains at most $2^n$ maximal finite subgroups up to conjugation.
\end{pr}

Let $F$ be a maximal finite subgroup in $\Gamma$.
Its fixed-point set $\mathrm{Fix}(F)$ is an affine subset of $\mathbb{E}^n$.
Furthermore, $\mathrm{Fix}(F)$ is mapped to a singular set $S_F$ in the quotient space $X=\mathbb{E}^n/\Gamma$.
Moreover, (1) a conjugation of $F$ does not change the set $S_F$; and, since $F$ is maximal, (2) $S_F$ is a \emph{simple} singular set;
that is, $S_F$ does not have a proper singular subset.
So instead of counting subgroups up to conjugation, we may count simple singular sets in $X$.
Crude bounds on this number were obtained by Tadashi Fujioka \cite{fujioka2025}.

It is expected that if the number of simple singular sets is maximal, then each of them has dimension 0;
that is, each set $S_F$ is an isolated point in $X$.
If this is indeed true, then the problem essentially boils down to the following well-known problem in discrete geometry:
\emph{show that one can choose at most $2^n$ points in $\mathbb{E}^n$ such that any triangle with vertices among the chosen points is acute or right.}
This problem was posed by Paul Erdős and solved by Ludwig Danzer and Branko Grünbaum \cite{danzer-gruenbaum}.
More on the subject can be found in the paper by Nina Lebedeva \cite{lebedeva}.

\subsection{Shortcut in a connected set}

\begin{pr}
Let $p$ and $q$ be a pair of points in a compact connected subset $K\subset \RR^n$.
Is it possible to connect $p$ to $q$ by a curve $\gamma$ such that the length of $\gamma\setminus K$ is arbitrarily small?
\end{pr}

Be aware that there are compact connected sets which contain no nontrivial paths; for example, the \emph{pseudoarc}.

The problem remains open for $n\ge 3$.
In the case $n=2$, an affirmative answer was obtained by Taras Banakh \cite{banakh}.
His argument is based on the following claim:

\textit{Suppose $K$ is a connected compact set in the plane and $p,q\in K$.
 Then, given $\eps>0$ there is a sequence of compact connected subsets $K_0,\dots, K_n$
 with pairs of points $p_i,q_i\in K_i$ for each $i$ such that
(1) $\diam K_i<\eps$ for each $i$,
(2) $p_0=p$, $q_n=q$, and
(3) if $\ell_i$ denotes the line segment $[q_{i-1}p_i]$, then
 \[\length \ell_1+\dots+\length \ell_n<\eps.\]}

Assume that the claim is proved.
Then we can prove the statement by applying it recursively to each $K_i$, choosing a much smaller value of $\eps$ at each step.
Indeed, consider the closure $L$ of the set of all points $p_i$ and $q_i$ for all iterations; it is a closed subset in $K$.
By adding to $L$ the line segments $\ell_i$ for all the iterations we obtain a curve connecting $p$ to $q$ that spends an arbitrarily small length outside of $L$.
Hence the 2-dimensional case follows.

Now let us prove the claim.
Consider a coordinate grid $\Gamma$ on the plane with step $\tfrac\eps2$.
Since the intersection $\Gamma\cap K$ is compact, it can be covered by a finite collection of closed line segments $\{I_1,\dots,I_n\}$ in the grid such that the
\[\sum_i\length(I_i\setminus K)<\eps.\]

Take a minimal collection of sets $\{K_i\}$ such that each set $K_i$ is a connected component of the intersection of $K$ with a closed square of the grid and the union
\[\biggl(\bigcup_i I_i\biggr)\cup\biggl(\bigcup_j K_j\biggr)\]
is a connected set containing $p$ and $q$.
Observe that the collection $\{K_i\}$ is finite and it satisfies the conditions in the claim, assuming that $K_i$ are ordered correctly.
(The choice of points $p_i$ and $q_i$ is straightforward;
the line segments $\ell_j$ can be chosen to lie in the union of $\{I_i\}$.)

\subsection{Involution of the 3-sphere}

\begin{pr}
Suppose that a closed geodesic $\gamma$ is the fixed-point set of an isometric involution of $(\mathbb{S}^3,g)$.
Assume that the sectional curvature of $g$ is at least~$1$.
Is it true that $\length\gamma\le2\cdot\pi$?
\end{pr}

A nonsharp universal upper bound on $\length\gamma$ follows from a
result of Tadashi Fujioka \cite{fujioka2025}.
The same question can be asked about spherical polyhedral metrics with curvature at least 1 in the sense of Alexandrov.

If, instead of an isometric involution, we have an $\mathbb{S}^1$-action, then the answer is ``yes''.
It follows since the quotient space $(\mathbb{S}^3,g)/\mathbb{S}^1$ is a convex disc with sectional curvature at least 1,
and $\gamma$ projects isometrically to its boundary.

\subsection{Faceting Fubini--Study}

A \index{polyhedral metric}\emph{polyhedral metric} on a compact manifold is a length metric such that, for some triangulation, each simplex is isometric to a simplex in Euclidean space.
Such a metric is said to be nonnegatively curved if the total angle around each codimension-two simplex does not exceed $2\cdot\pi$; this is equivalent to nonnegative curvature in the sense of Alexandrov.

\begin{pr}
Let $d$ be the Fubini--Study metric on the complex projective plane $\CP^2$.
Is it possible to find a sequence of nonnegatively curved polyhedral metrics $d_n$ on $\CP^2$ that converges to $d$?
\end{pr}

Nonnegatively curved polyhedral metrics on $\CP^2$ have a number of interesting properties.
For example, any such metric is Kähler; that is, one can choose a complex structure on each simplex of its triangulation that is respected by the gluings.
Furthermore, these complex structures coincide with the standard complex structure on $\CP^2$ for a right parametrization, and in this parametrization, the singular set forms an algebraic submanifold.
Such metrics on $\CP^2$ were constructed and studied by Dmitri Panov \cite{panov-Kaeler};
see also his paper with Martin de Borbon \cite{deborbon-panov2026}.

\subsection{Octahedron comparison}

Suppose $\Gamma$ is a graph with vertices $v_1,\dots,v_n$.
We say that a metric space $X$ meets the \emph{$\Gamma$-comparison} if, for any set of points in $X$ labeled by the vertices of $\Gamma$, there is a model configuration $\tilde v_1,\dots,\tilde v_n$ in the Hilbert space $\HH$ such that
if $v_i$ is adjacent to $v_j$, then
$|\tilde v_i-\tilde v_j|_{\HH}\z\le | v_i-v_j|_{X}$,
and otherwise,
$|\tilde v_i-\tilde v_j|_{\HH}\ge | v_i-v_j|_{X}$;
here $|p-q|_M$ denotes the distance between points $p$ and $q$ in the metric space~$M$.


{

\begin{wrapfigure}{r}{30 mm}
\vskip-6mm\centering
\begin{tikzpicture}[x=1cm,y=1cm]
\coordinate (z0) at ( 0,       0.5);\coordinate (z1) at (-0.7794,  0.45);\coordinate (z2) at (-0.4330, -0.25);\coordinate (z3) at ( 0,      -0.9);\coordinate (z4) at ( 0.4330, -0.25);\coordinate (z5) at ( 0.7794,  0.45);
\draw[use Hobby shortcut, closed=true]  (z0) .. (z1) .. (z3) .. (z4);\draw[use Hobby shortcut, closed=true]  (z1) .. (z2) .. (z4) .. (z5);\draw[use Hobby shortcut, closed=true]  (z2) .. (z3) .. (z5) .. (z0);
\fill (z0) circle (1.2pt);\fill (z1) circle (1.2pt);\fill (z2) circle (1.2pt);\fill (z3) circle (1.2pt);\fill (z4) circle (1.2pt);\fill (z5) circle (1.2pt);
\node[above]  at (z0) {$x$};\node[above,xshift=1pt] at (z3) {$x'$};
\node[below right] at (z4) {$y$};\node[above left,xshift=3pt] at (z1) {$y'$};
\node[below left]  at (z2) {$z$};\node[above right] at (z5) {$z'$};
\end{tikzpicture}
\end{wrapfigure}

Let us denote by $O_3$ the graph of the octahedron.
For the labeling shown on the diagram, $O_3$-comparison means that any six-point configuration $x,x',y,y',z,z'\in X$ can be mapped to $\HH$ so that the diagonals $xx'$, $yy'$, and $zz'$ do not get shorter and the edges do not get longer.

}

\begin{pr}
Is it true that $O_3$-comparison holds in any geodesic $\CAT(0)$ space?
\end{pr}

This problem is motivated by the question of Michael Gromov \cite[1.19$_+(e)$]{gromov-MetStr}.

Note that the 4-cycle $C_4$ is an induced subgraph of $O_3$.
Therefore, \textit{$O_3$-comparison implies $C_4$-comparison};
the latter describes the $\CAT(0)$-comparison \cite{akp}.
In general, $C_4$-comparison does not imply $O_3$-comparison,
but for intrinsic metrics it is unknown; a weaker inequality was proved by Tetsu Toyoda \cite{toyoda2025}.

An affirmative answer might lead to a classification of six-point spaces that admit an isometric embedding in a geodesic $\CAT(0)$ space.
The five-point case was solved by Tetsu Toyoda \cite{toyoda2020} (see also \cite{lebedeva-petrunin2021}).

On the other hand, Nina Lebedeva and the author proved that $O_3$-comparison holds in products of trees \cite{lebedeva-petrunin2025}.
Therefore, a counterexample would mean that $O_3$-comparison is something new and interesting;
so both answers are good.

\parbf{Postscript.}
The problem has been solved by Jihun Kim \cite{kim2026}.

\section*{Overheard}

\subsection{Group generated by central symmetries}
\label{Group generated by central symmetries}

\begin{pr}
Consider the Lobachevsky space of sufficiently large dimension.
Does it admit a discrete cocompact isometric action generated by
\begin{enumerate}[(a)]
 \item elements of finite order?
 \item central symmetries?
 \item\label{vinberg} rotations around subspaces of codimension $k$?
\end{enumerate}

\end{pr}

Ernest Vinberg \cite{vinberg, vinberg-strong} proved that there are no such actions generated by reflections.
This solves part (\ref{vinberg}) for $k=1$; the remaining questions are completely open.

\subsection{Milnor's cartography problem}
\label{Milnor's cartography problem}

\begin{pr}
Let $\Omega$ be a round disc of radius $\alpha<\tfrac\pi2$ on the unit sphere $\mathbb{S}^2$,
and let $\Omega'$ be a convex domain in $\mathbb{S}^2$ with the same area as $\Omega$. Show that there is a $(1,\tfrac{\alpha}{\sin\alpha})$-bi-Lipschitz map from $\Omega'$ to the plane.
\end{pr}

The question was asked by John Milnor \cite{milnor-cartography}.
Polar coordinates with the origin at the center of $\Omega$ provide an example of a $(1,\tfrac{\alpha}{\sin\alpha})$-bi-Lipschitz map and this constant $\tfrac{\alpha}{\sin\alpha}$ is optimal.

If, in the formulation of the problem, one replaces ``area''  with ``perimeter'',
then the answer is ``yes''.
It can be proved using the idea of Gary Lawlor \cite{lawlor}.

\parbf{Postscript.}
The problem has since been solved by the author \cite{petrunin2026}.

\subsection{Convex hull in CAT(0)}
\label{Convex hull in CAT(0)}

Let $X$ be a geodesic metric space.
Recall that a set $C\subset X$ is called \emph{convex} if, for any two points $x,y\in C$, every geodesic from $x$ to $y$ lies in~$C$.
The \emph{closed convex hull} of a subset $K\subset X$ is defined as the intersection of all closed convex sets containing $K$.

\begin{pr}
Let $K$ be a compact subset of a complete length $\mathop{\rm CAT}(0)$ space $X$.
Is it true that the closed convex hull of $K$ is compact?
\end{pr}

Note that the space is not assumed to be locally compact; otherwise, the statement is evident.
I believe in a counterexample, perhaps even with negatively pinched curvature (in the sense of Alexandrov).

The closed convex hull of a compact set (in particular, of a three-point set) might look ugly even in a smooth 3-dimensional manifold.
If in doubt, check the note of Alexander Lytchak and the author \cite{lytchak-petrunin}.

This problem was mentioned by Mikhael Gromov \cite[6.B$_1\text{(f)}$]{gromov-asymptotic}
and restated by Eva Kopeck\'a and Simeon Reich \cite{kopecka-reich}, 
but it seems to have been known in mathematical folklore for a long time.

The problem remains open even when $K$ consists of three points.
A counterexample is known if the $\mathop{\rm CAT}(0)$ assumption is weakened to the existence of a conical geodesic bicombing; it was constructed by Giuliano Basso, Yannick Krifka, and Elefterios Soultanis \cite{BKS}.

\subsection{Homeomorphism of the torus}

Let $\TT$ denote the 2-dimensional torus.

A homeomorphism $h\:\TT\to \TT$ is called \emph{simple} if it is supported on an embedded open disc;
that is, there is an embedded open disc $\Delta\subset \TT$ such that $h$ is the identity outside of $\Delta$.

The \emph{disc norm} of a homeomorphism $h\:\TT\to \TT$ is defined as the least number of simple homeomorphisms $h_1,\dots,h_n$ such that
\[h=h_1\circ\dots\circ h_n.\]

\begin{pr}
Is there a homeomorphism of the 2-dimensional torus with arbitrarily large disc norm?
\end{pr}

Analogous questions for the cylinder and the sphere are solved, but for the Möbius band the problem is open.
This problem is discussed by Dmitri Burago, Sergei Ivanov, and Leonid Polterovich in \cite{BIP}.

\parbf{Postscript.}
Jonathan Bowden, Sebastian Hensel, and Richard Webb \cite{bowden-hensel-webb} constructed unbounded quasimorphisms on the group of homeomorphisms of the torus that are isotopic to the identity, which implies an affirmative answer to the problem.

\subsection{Open projection}

Let us denote by $\DD^n$ a closed $n$-dimensional ball.

Recall that a map is called \emph{open} if the image of any open set is open.

\begin{pr}
Is it possible to construct an embedding $\DD^4\hookrightarrow \mathbb{S}^2\times
\mathbb R^2$
such that the projection  $\DD^4\to \mathbb{S}^2$ is an open map?
\end{pr}

It is easy to construct an embedding $\DD^3\hookrightarrow \mathbb{S}^3$
such that its composition with the Hopf fibration is open.
Therefore, there is an open map $f_3\:\DD^3\to\mathbb{S}^2$.

Composing $f_3$ with any open map $\DD^4\to \DD^3$,
one gets an open map $f_4:\DD^4\to \mathbb{S}^2$.
However, this map is not a projection of an embedding $\DD^4\hookrightarrow \mathbb{S}^2\times
\mathbb R^2$.

The question was asked by $\eps$-$\delta$ (an anonymous mathematician)  \cite{epsilon-delta}.

\subsection{Zero-curvature along geodesics}

\begin{pr}
Let $M$ be a surface of genus 2 with a Riemannian metric of nonpositive curvature.
Show that, with probability 1, a random geodesic in $M$ runs into the set of negative curvature.
\end{pr}

I learned this problem from Dmitri Burago; it originated from a discussion with Keith Burns, Michael Brin, and Yakov Pesin \cite[see also][]{hertz}.
Weisheng Wu \cite{wu2015} proved the statement under the additional assumption that the set of points with negative curvature has finitely many connected components.
It seems that \textit{there are no known examples of a surface $M$ with a nonperiodic geodesic that runs in the set of zero curvature}.

\subsection{Ambrose problem}

\begin{pr}
Let $(M,g)$ be a simply connected complete Riemannian manifold, and let $p\in M$.
Denote by $h$ the pullback of the metric tensor $g$ to the tangent space $\T_pM$ along the exponential map $\exp_p\:\T_pM\z\to M$;
in general, $h$ might degenerate at some points.
Does $h$ define the manifold $(M,g)$ and the point $p$ up to isometry? 
\end{pr}

This simple-looking question was asked a long time ago by Warren Ambrose \cite{ambrose}.
Affirmative answers were obtained in three cases: (1) in dimension 2, (2) for analytic metrics, and (3) for generic metrics; see \cite{hebda1994, hebda2010, itoh, ardoy} and the references therein.

\subsection{Optical fibers}

Suppose that $V$ is a body in $\RR^3$ bounded by two plane figures $F_1$ and $F_2$ (the ends) and a smooth surface $\Sigma$ (a tube).
Suppose that a light ray entering $V$ thru one end and bouncing with perfect reflection from the interior walls of $\Sigma$ will emerge from the other end with probability 1;
moreover, a light ray that starts in the interior of $V$, after bouncing, emerges from one of the ends with  probability 1.
In this case, we say that $V$ is an \emph{optical fiber}.

A small tubular neighborhood of a smooth curve $\gamma\:[a,b]\to \RR^3$ with disk ends provides an example of an optical fiber.
Let us describe a more general example, which can be constructed in the same way using any plane figure bounded by a smooth curve instead of a disk.
Suppose a plane curve is given by a periodic parametrization $\theta\mapsto (x(\theta),y(\theta))$.
Choose a  parallel normal frame $e_1,e_2$ along $\gamma$
(that is, $e_i'(t)\parallel\dot\gamma(t)$ for all~$t$),
and consider the tube $[a,b]\times\mathbb S^1\hookrightarrow\mathbb R^3$ defined by
$$(t,\theta)\mapsto \gamma(t)+x(\theta){\cdot}e_1(t)+y(\theta){\cdot}e_2(t)$$
(The condition that the frame is parallel implies that any normal plane to $\gamma$ cuts the tube at a right angle.)
This way we get an \emph{optical fiber} with congruent ends.

Are there any constructions of optical fibers different from the one described above? In other words:

\begin{pr}
Is it always possible to slice an optical fiber into plane figures that meet the walls at right angles?
\end{pr}

In particular, 

\begin{pr}
Is there an optical fiber with noncongruent ends?
\end{pr}

By  Liouville's theorem, it is clear that the ends must have the same area.

If the walls are only assumed to be \emph{piecewise} smooth, then one can construct an optical fiber with a pair of equidecomposable figures at the ends.
(The construction is the same, but one splits the tube into a few tubes along the way and then reassembles them.)

The problem is based on two questions posted on MathOverflow \cite{rourke-optic, petrunin-optic}.

\subsection{Guth's sponge}

\begin{pr}
Show that there is $\eps_n>0$ such that any bounded open set in $\RR^n$ with volume at most $\eps_n$ admits a length-increasing embedding into a unit ball.
\end{pr}

The question was asked by Larry Guth in his thesis \cite{guth-2005}.
It is open in all dimensions $n\ge 2$.
In dimension 2, Parker Glynn-Adey proved a weaker version using a strip between parallel lines instead of the ball \cite{glynn-adey}.



\sloppy
\printbibliography[heading=bibintoc]
\fussy

\end{document}